# Optimal Operation of PV-Diesel MicroGrid with Multiple Diesel Generators under Grid Blackouts

Mansour Alramlawi, *Student Member, IEEE*, Ayeni Femi Timothy, Aouss Gabash, *Member, IEEE*,
Erfan Mohagheghi, *Member, IEEE*, Pu Li
Department of Simulation and Optimal Processes
Institute of Automation and Systems Engineering
Ilmenau University of technology
Ilmenau, Germany
mansour.alramlawi@tu-ilmenau.de, femi-timothy.ayeni@tu-ilmenau.de, aouss.gabash@tu-ilmenau.de,
erfan.mohagheghi@tu-ilmenau.de, pu.li@tu-ilmenau.de

*Abstract*—This paper addresses the optimal operation problem of a PV-diesel microgrid considering grid blackouts, which is a usual case of discontinuous power supply in developing countries. The model of a grid-connected PV-diesel microgrid is enhanced, and new practical constraints are added. In addition, a new mixed-integer nonlinear programming (MINLP) problem is formulated to optimize the power flow in the microgrid for covering the load while minimizing the power consumption cost as well as maximizing the dispatched power from the PV-array. To demonstrate the applicability of the developed model and the optimal operation strategy, a comparison between different system configurations is made. The results show that the proposed microgrid can reduce the dispatched power cost while decreasing the power curtailment from the PV-array. Moreover, installing multiple diesel generators (DGs) can significantly reduce the fuel consumption by the DGs in comparison to a single DG in the PV-diesel microgrid.

## I. INTRODUCTION

Diesel generators are widely used as backup power supplies to cover the loads during grid blackouts [1], which is a major problem in many countries worldwide [2], [3]. Due to the technology development in the recent years, it is possible to integrate renewable sources with conventional sources to build a local microgrid [4], [5]. However, many technical and economical complexities have to be handled in a microgrid operation to provide reliable and cost-effective energy.

The optimal design and operation of a PV-diesel-battery microgrid have been heavily investigated in the literature [6], [7]. However, battery banks have some drawbacks such as high capital cost, replacement cost, and battery efficiency reduction. Therefore, a battery-less PV-diesel microgrid is used as an economical solution with a low capital cost to provide uninterrupted power supply to the load, while reducing the fuel consumption of a diesel generator. An experimental study was done in [8] to measure the performance of a

This work is supported by the German Academic Exchange Service (DAAD).

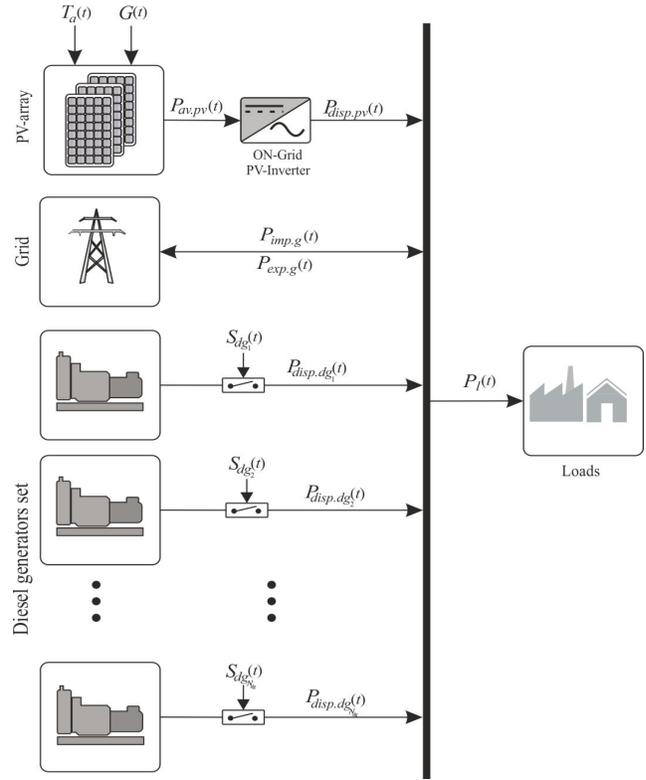

Fig. 1. Schematic diagram of the proposed hybrid PV-diesel system

PV-diesel system without battery under different operating conditions. The results show that the system operation should be optimized to enhance the system efficiency. In another work [9], the performance of a battery-less PV-diesel microgrid was evaluated over one year period, and recommendations were given to increase the system profitability. In [10], the size of a battery-less PV-diesel system was optimized to decrease the net capital cost of the system and the consumed energy cost. Moreover, the work in [11] presented an optimization approach to find the optimal surface area of the PV-array and the optimal capacity of the diesel generator using a harmony search algorithm. A flexible operation strategy was proposed in [12] to track the maximum power from the PV-array, regulate load voltage and control the diesel generator engine in a PV-diesel system but without optimizing the power flow in the system. From another perspective, the work in [13] proposed an automatic control strategy for reactive power generation in an isolated wind-diesel system.

A key limitation of the previous studies is that the optimal operation problem of the grid-connected PV-diesel microgrid with multiple diesel generators under grid blackouts was not investigated. Moreover, a special property of the grid feeding PV-inverters is that they cannot operate in the island mode if there is no grid forming device for setting the voltage and frequency in the microgrid [14], [15]. This leads to a new constraint that should be held in dispatching the generated power from the PV-array. In this paper, the model of the grid-connected battery-less PV-diesel microgrid is enhanced and new practical constraints are added. In addition, a new approach is developed to optimize the operation of the PV-diesel microgrid which, simultaneously minimizes the energy consumption costs and maximizes the dispatched power from the PV-array.

The remainder of the paper is organized as follows. In Section II, the modeling of the proposed PV-diesel microgrid is explained. Section III describes the optimization problem. In Section IV, the results of the paper are discussed . The paper is concluded in Section V.

## II. PV-DIESEL MICROGRID DESCRIPTION

Fig.1 illustrates the considered microgrid in this work. It is a grid connected microgrid consisting of a PV-array, a grid feeding PV-inverter, multiple diesel generators, controllable switches and loads. The mathematical model of each part is explained and given below.

### A. PV-Array

In this work, the three component model of the solar irradiance [16] is used to calculate the total solar irradiance $G_t(t)$ arriving at an inclined PV-module.

$$G_t(t) = G_B(t) + G_D(t) + G_R(t). \quad (1)$$

Here $G_B(t)$ is the direct beam irradiance that arrives the PV-module without reflection or scattering, $G_D(t)$ is the diffused solar irradiance which is scattered by the clouds and $G_R(t)$ is the reflected solar irradiance by the ground, see Fig. 2. The mentioned solar irradiance components are calculated as follows

$$G_B(t) = R_B G_b(t) \quad (2)$$

$$G_D(t) = \frac{(1+cos\beta)}{2} G_d(t) \quad (3)$$

$$G_R(t) = \frac{(1-cos\beta)}{2} \rho G(t) \quad (4)$$

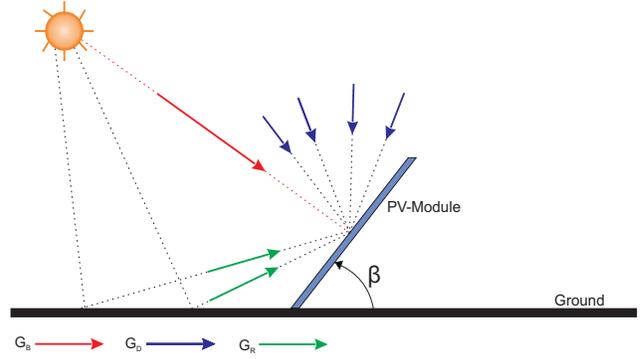

Fig. 2. Solar irradiance components that arriving at an inclined PV-module. [16]

where $G(t)$ is the global solar irradiance, $G_b(t)$ and $G_d(t)$ are the direct beam and diffused solar irradiance measured at the horizontal surface. $R_B$ is the direct beam irradiance factor, $\beta$ is the PV-module inclination angle and $\rho$ is the reflection factor that depends on the ground surface type.

Additionally, the mathematical model of a single diode PV-cell with a series resistance is used to calculate the maximum output power from a PV-cell $P_{max.c}(t)$, which depends on the instantaneous value of the solar irradiance and the ambient temperature, as follows [3]

$$P_{max.c}(t) = V_{oc.c}(t) I_{sc.c}(t) FF(t) \quad (5)$$

$$V_{oc.c}(t) = V_{oc.c.stc} + K_v(T_c(t) - 25) \quad (6)$$

$$I_{sc.c}(t) = (I_{sc.c.stc} + K_i(T_c(t) - 25)) \frac{G(t)}{1000} \quad (7)$$

$$FF(t) = \frac{P_{max.c}(t)}{V_{oc.c}(t) I_{sc.c}(t)} \quad (8)$$

where $V_{oc.c}(t)$ and $I_{sc.c}(t)$ are the open circuit voltage and the short circuit current of the PV-cell, receptively. $V_{oc.stc.c}(t)$ and $I_{sc.stc.c}(t)$ are the open circuit voltage and the short circuit current of the PV-cell at standard test conditions, receptively. $K_v$ is the open circuit voltage temperature coefficient, $Ki$ is the short circuit current temperature coefficient, $T_c$ is the PV-cell temperature coefficient and $FF(t)$ is the PV-cell fill factor. The total available power from the PV-array $P_{av.pv}(t)$ is calculated by

$$P_{av.pv}(t) = N_{s.m} N_{p.m} N_{c.m} P_{max.c}(t) \quad (9)$$

where $N_{s.m}$ and $N_{p.m}$ are the number of series and parallel connected modules in the PV-array. $N_{c.m}$ is the number of PV-cells in PV-array.

## B. Diesel Generator

The fuel consumption of each diesel generator $f_{con.dg_i}(t)$ is related to the power dispatched from it $P_{disp.dg_i}(t)$ and its rated power $P_{r.dg_i}$ as follows [17]

$$f_{con.dg_i}(t) = \begin{cases} AP_{disp.dg_i}(t) + BP_{r.dg_i}, & if\ S_{dg_i}(t)=1 \\ 0 & otherwise \end{cases} \quad (10)$$

where $A$ is a constant with the value of $0.246 l/kWh$, $B$ is a constant with the value $0.08415 l/kW$. The start-up and shutdown costs of each diesel generator are calculated as follows

$$\alpha_{up.dg_i}(t) = S_{dg_i}(t)(1 - S_{dg_i}(t-1)),\ \forall i=1,2,...,N_{dg} \quad (11)$$

$$\alpha_{d.dg_i}(t) = S_{dg_i}(t-1)(1 - S_{dg_i}(t)),\ \forall i=1,2,...,N_{dg}. \quad (12)$$

Here $\alpha_{up.dg_i}(t)$ and $\alpha_{d.dg_i}(t)$ are auxiliary binary variables that represent the changes at each diesel generator generator status. If $\alpha_{up.dg_i}(t) = 1$, it means that the $i^{th}$ diesel generator is started up at this time step, otherwise change is occurred. If $\alpha_{d.dg_i}(t) = 1$, it means that the $i^{th}$ diesel generator is shut down at this time step itherwise no change is occurred. The frequent operation of the diesel generator at law load increase the risk of diesel generator engine failure due to wet stacking problem [9]. Therefore, the following constraint is added

$$P_{disp.dg_i}(t) \geq 0.3 P_{r.dg_i}, \qquad \forall i=1,2,...,N_{dg} \quad (13)$$

In this work, it is assumed that a grid feeding PV-inverter is used to inject the generated power from the PV-array in the microgrid. It is known that the grid feeding PV-inverters are unable to work in an island mode in the case of a blackout [14], [15]. Therefore, at least one of the diesel generators should be activated to set the voltage and frequency of the microgrid. For this, the following constraint should be held during the blackouts periods

$$\sum_{i=1}^{N_{dg}} S_{dg_i}(t) \geq 1, \qquad \forall i=1,2,...,N_{dg}. \quad (14)$$

Here $S_{dg_i}(t)$ is the status of each diesel generator switch.

## C. Grid Blackouts

In this paper, it is assumed that the proposed microgrid is connected to a public grid suffering from scheduled blackouts [3], [18]. The grid capability $P_g(t)$ is described as

$$P_g(t) = \alpha_g(t) P_{g.max} \quad (15)$$

where $P_{g.max}$ is the grid maximum allowed power to be imported or exported from/to the grid. $\alpha_g(t)$ represents the status of the grid, i.e. when $\alpha_g(t) = 1$, the grid is ON and when $\alpha_g(t) = 0$, the grid is OFF. The behavior of $\alpha_g(t)$ is illustrated in Fig. 3, where $\tau$ is the Grid-ON period and $T$ is the total ON-OFF period.

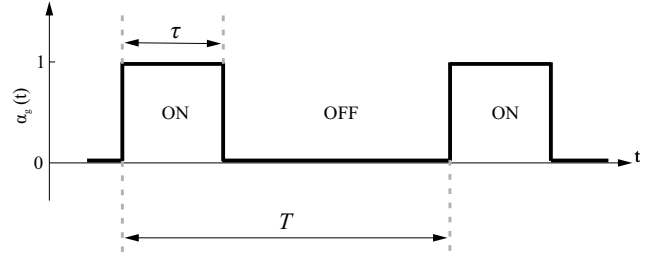

Fig. 3. The illuration of grid status [3].

Moreover, it is assumed that the microgrid is able to import and export active power to the public grid [19]–[21]. Meanwhile, the following constraints should be held

$$\delta_{imp}(t) P_{disp.g}(t) \leq P_g(t) \quad (16)$$

$$\delta_{exp}(t) P_{exp.g}(t) \leq P_g(t) \quad (17)$$

where

$$\delta_{imp}(t) + \delta_{exp}(t) \leq 1 \quad (18)$$

where $\delta_{imp}(t)$ and $\delta_{exp}(t)$ are binary variables used to prevent importing and exporting active power from/to the public grid simultaneously.

## III. OPTIMIZATION PROBLEM FORMULATION

The objective of the proposed optimal operation strategy is to minimize the cost of the consumed energy by the load and meanwhile to maximize the dispatched power from the PV-array. The objective function is defined as follows

$$\min_{\mathbf{u}_c(t),\mathbf{u}_b(t)} J(t) = w_1(C_{e.g} P_{disp.g}(t) + \sum_{i=1}^{N_{dg}}(C_f f_{con.dg_i}(t) \\ + C_{up}\alpha_{up.dg_i}(t) + C_d \alpha_{d.dg_i}(t) + C_{o\&m} S_{dg_i}(t))) \\ - w_2(P_{disp.pv}(t) + C_{exp} P_{exp}(t)). \quad (19)$$

Here, $\mathbf{u}_c(t)$ and $\mathbf{u}_b(t)$ are the continuous and binary control variables vectors, respectively. $w_i$ is a weighting factor used to convert the multi-objective optimization problem to a single-objective optimization problem [22]. $C_{e.g}$ is the cost of energy dispatched from the grid in $/kWh$, $C_{exp}$ is the cost of the exported power to the grid in $/kWh$, $C_f$ is the cost of diesel generator fuel in $/l$, $C_{up/d}$ is the start up or shut down cost in $ and $C_{o\&m}$ is the maintenance cost in $/h$. $P_{disp.pv}(t)$ is the dispatched power from the PV-array to cover the load, $P_{exp}(t)$ is the exported power from the PV-array to the grid.

In addition to the new constraints in Eq.14 and Eqs.16-18 the operational constraints of the PV-array, battery and the grid stated in [3] should be held.

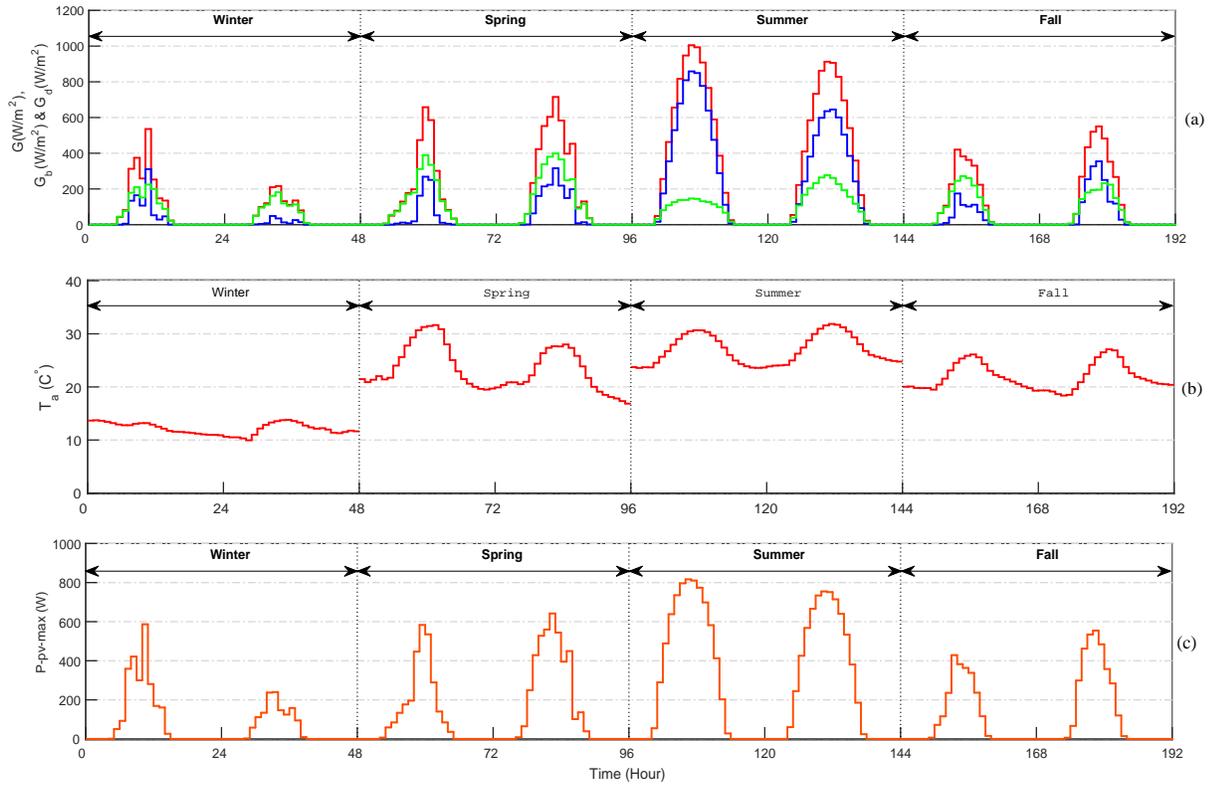

Fig. 4. PV-array model illustration. (a) Solar irradiance components at horizontal surface : global irradiance (red), direct beam irradiance (green) and diffused irradiance (blue). (b) Ambient temperature. (c) The maximum generated power by the PV-module.

## IV. RESULTS AND DISCUSSIONS

To illustrate the potential of the proposed microgrid, a case study is adapted from [3]. To demonstrate the benefits of using a PV-diesel microgrid with multiple diesel generators, three case studies are investigated. In case 1, a stand-alone diesel generator with size 500 $kVA$ is used to cover the load during blackouts periods. In case 2, a 700 $kW_p$ PV-array is added to the system. In case 3, three diesel generators with 200, 150 and 150 $kVA$ are used instead of one diesel generator. The solar irradiance components and the ambient temperature for selected two days from each season are shown in Fig. 4 (a) and (b), respectively. The maximum generated power from 1 $kW_p$ PV-array is shown in Fig. 4 (c). The load profile is addapted from [23] and shown in Fig. 5 (a). The computation results for one year are shown in Tables I and II.

It is clearly seen from Fig. 5 that the proposed microgrid is able to cover the load while satisfying the microgrid technical constraints. Due to considering the cost of the dispatched energy from the grid and the diesel generator in the optimization problem, a huge reduction of the total dispatched energy from the grid $E_{disp.g}$ and the diesel generator $E_{disp.dg}$ is gained, as given in Table I. Furthermore, the generated power from the PV-array is used to cover the load and the extra power is exported to the grid as shown in Fig. 5 (c). From Fig. 5

TABLE I
YEARLY DISPATCHED ENERGY FROM THE GRID AND DIESEL GENERATOR AT EACH CASE.

|  | Case 1 | Case 2 | Case 3 |
| --- | --- | --- | --- |
| $E_{disp.g}(MWh)$ | 1177.3 | 667.32 | 667.32 |
| $E_{disp.dg}(MWh)$ | 1176.9 | 919.98 | 770.59 |
| $E_{disp.pv}(MWh)$ | N/A | 766.18 | 915.60 |

TABLE II
COST OF THE YEARLY DISPATCHED ENERGY AT EACH CASE.

|  | Case 1 | Case 2 | Case 3 |
| --- | --- | --- | --- |
| $C_{disp.g}(k\$)$ | 176.59 | 100.1 | 100.1 |
| $C_{disp.dg}(k\$)$ | 685.66 | 586.45 | 382.98 |
| Total $(k\$)$ | 862.25 | 686.55 | 483.08 |

(d), it can be seen that during blackout periods at least one of the diesel generators is ON to satisfy the constraint in Eq. 14. Meanwhile, the diesel generator operation status is changed based on the load and the generated power from the PV-array level, to decrease the diesel generators fuel consumption and to increase the total dispatched energy from the PV-array. The allowed power to be dispatched from the PV-array is limited to satisfy the minimum dispatched power from the diesel generator constraint, see Eq. 13. Therefore, using multiple diesel generators with different sizes increase the allowed dispatched power from the PV-array be selecting the optimal diesel generator to be turned ON. This decreases the dispatched energy from the diesel generators as shown in Table I. Moreover, the cost of the dispatched energy from the grid and the diesel generator is given in table II. It can be seen that using the PV-diesel microgrid with multiple diesel generators (case 3), reduces the power consumption cost up to 29.64% in comparison to case 2 and 43.97% in comparison to case 1.

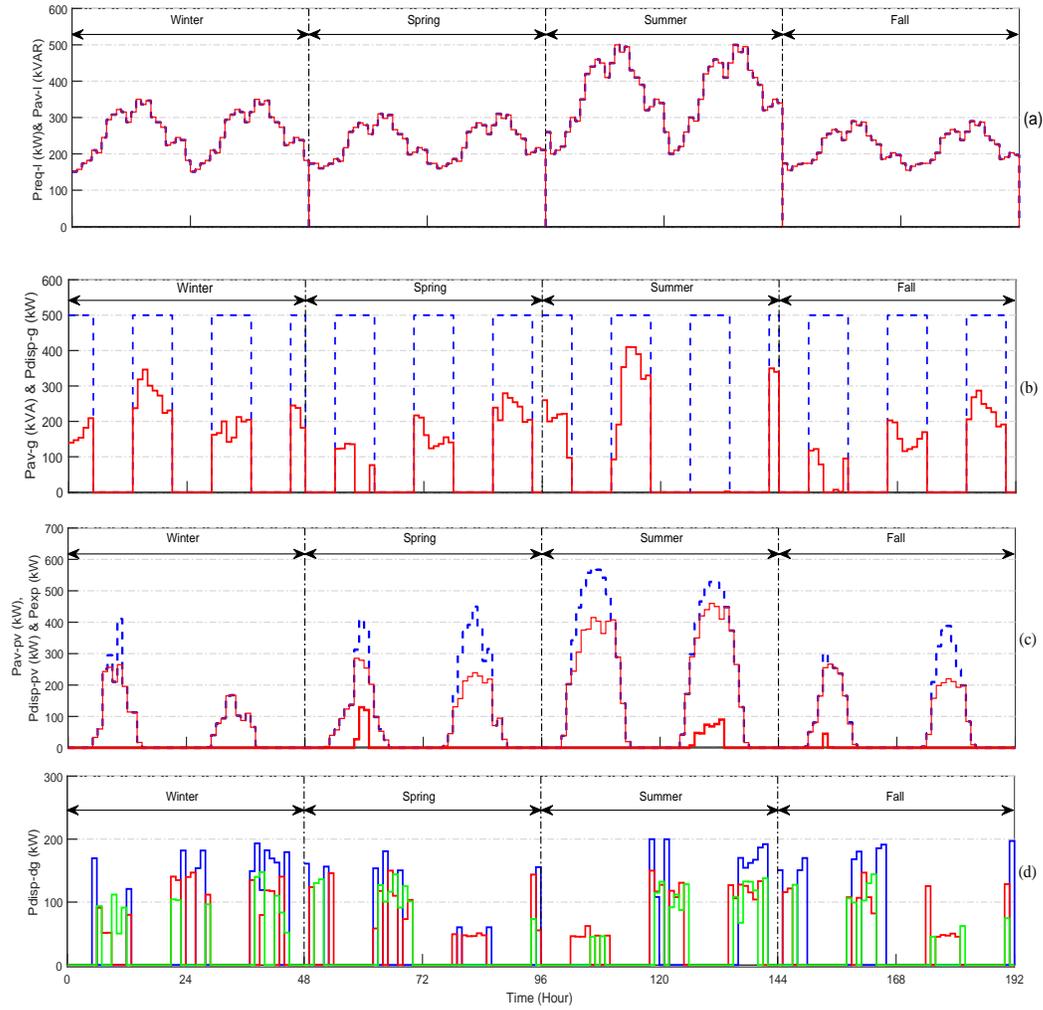

Fig. 5. Microgrid operation in arbitrary two days for each season in case 3. (a) Required load power (dashed-blue) and the available power (solid-red). (b) Available power (dashed-blue) and dispatched power (solid-red) from grid. (c) Available power (dashed-blue), dispatched power to the load (thin-red) and exported power (bold-red) from PV-array. (d) Dispatched power from $200kVA$ diesel generator (blue), dispatched power from the first $150kVA$ diesel generator (green) and dispatched power from the second $150kVA$ diesel generator (red).

## V. CONCLUSION

In this work, a battery-less PV-diesel microgrid is developed to cover the loads during blackouts periods instead of using only diesel system. The results show that the proposed microgrid can provide an uninterrupted power supply to the load and considerably reduce the total dispatched energy cost. In addition, an optimal operation strategy is introduced to reduce the power dispatched from the grid and the diesel generator while decreasing the power curtailment from the PV-array, which increases the economic and environmental benefits of installing the PV-array. Moreover, it is proved that using multiple diesel generators a battery-less PV-diesel microgrid is more efficient than using single diesel generator. On the basis of the promising findings presented in this paper, future work will involve the combined optimal design and operation problem of the developed microgrid.